\documentclass[10pt]{amsart}
\usepackage{amssymb}
\theoremstyle{plain}
\newtheorem{thm}{Theorem}
\theoremstyle{definition}
\newtheorem{prop}{Proposition}[section]
\newtheorem{cor}[prop]{Corollary}
\newtheorem{lem}[prop]{Lemma}

\newtheorem{defn}[prop]{Definition}
\newtheorem{claim}[prop]{Claim}
\theoremstyle{remark}

\newcommand{\D}{\mathfrak{D}}

\DeclareMathOperator{\ord}{Ord}
\DeclareMathOperator{\dom}{dom}

\DeclareMathOperator{\ran}{ran}

\DeclareMathOperator{\cf}{cf}

\newcommand{\sk}{\vskip.05in}

\newcommand{\restr}{\upharpoonright}
\newcommand{\forces}{\Vdash}

\newcommand{\subs}{\subseteq}

\numberwithin{equation}{section}

\begin{document}
\title{On iterated forcing for successors of regular cardinals}
\author{Todd Eisworth}
\address{Department of Mathematics\\
         University of Northern Iowa\\
         Cedar Falls, IA\\
         50614}
\email{eisworth@math.uni.edu}
\begin{abstract}
We investigate the problem of when $\leq\lambda$--support iterations of $<\lambda$--complete notions of forcing preserve $\lambda^+$.  We isolate a property --- {\em properness over diamonds} --- that implies $\lambda^+$ is preserved and show that this property is preserved by $\lambda$--support iterations. Our condition is a relative of that presented in \cite{RoSh:655}; it is not clear if the two conditions are equivalent. We close with an application of our technology by presenting a consistency result on uniformizing colorings of ladder systems on $\{\delta<\lambda^+:\cf(\delta)=\lambda\}$ that complements a  theorem of Shelah in \cite{Sh:f}. 
\end{abstract}
\keywords{}
\subjclass{}
\date{\today}
\thanks{The author would like to thank the Graduate College of the University of Northern Iowa for a Summer Fellowship in 2002 that provided support while this research was carried out. The author would also like to thank Andrzej Roslanowski for listening to preliminary versions of the proofs in this paper.}
\maketitle
\section{Definitions}

One of the mysteries of iterated forcing theory is the lack of a good solution to the following ``equation'' for an uncountable regular cardinal $\lambda$:
\begin{equation*}
\frac{\text{proper forcing}}{\text{countable support iteration}}=\frac{x}{\text{$\lambda$--support iteration}}.
\end{equation*}
The goal of this paper is to present a generalization of properness to the context of larger cardinals. We make no claim that ours is the ``right'' generalization; however, the proof that our condition is preserved by $\lambda$--support iteration is close to the proof that properness is preserved by countable support iteration and seems quite natural.

Throughout this paper, we make the following assumptions:

\begin{itemize}
\item $\lambda$ is a regular cardinal satisfying $\lambda=\lambda^{<\lambda}$.
\item $\D$ is a normal filter on $\lambda$ ``with diamonds'', i.e., for every $S\in \D^+$, there is a sequence $\langle A_\delta:\delta\in S\rangle$ such that for every $A\subs\lambda$,
\begin{equation*}
\{\delta\in S: A\cap \delta= A_\delta\}\in \D^+.
\end{equation*}
\item $\chi$ is a regular cardinal that is ``large enough''.
\end{itemize}

We are going to be looking at when $\lambda^+$ is preserved by $(\leq)\lambda$--support iterations of $(<)\lambda$--complete notions of forcing. Just as in the case of proper forcing, we will have to look at how our forcing notions interact with elementary submodels.

\begin{defn}
Let $N$ be an elementary submodel of $H(\chi)$.  We say that $N$ is {\em relevant} if 
\begin{itemize}
\item $||N||=\lambda$
\item $N^{<\lambda}\subs N$
\item $N=\bigcup_{\alpha<\lambda}N_\alpha$, where $\langle N_\alpha:\alpha<\lambda\rangle$ is a continuous $\in$--increasing sequence of elementary submodels of $H(\chi)$ such that $\langle N_\beta:\beta\leq\alpha\rangle\in N_{\alpha+1}$ and $|N_\alpha|<\lambda$. (We say that $\langle N_\alpha:\alpha<\lambda\rangle$ is a {\em filtration} of $N$.)
\end{itemize}
\end{defn}

The natural attempt at generalizing properness results in a definition along the following lines:
\begin{defn}
\label{properdef}
A notion of forcing $P$ is said to be $\lambda$--proper if for all sufficiently large regular cardinals ~$\chi$, there is some $x\in H(\chi)$ such that whenever $M$ is a relevant elementary submodel of  $H(\chi)$ with $\{P, x\}\in M$ and
$p$ is an element of $M\cap P$, there is a condition $q\leq p$ such that
\begin{equation*}
q\forces\text{``} M[\dot G_P]\cap\ord = M\cap\ord\text{''}.
\end{equation*}
Such a condition $q$ is said to be $(M, P)$--generic.
\end{defn} 

Some of the qualities of properness generalize in a straightforward fashion to this new context. For example,  $\lambda$--proper notions of forcing do not collapse $\lambda^+$, and it is easy to prove that both $\lambda^+$--closed and $\lambda^+$--c.c. notions of forcing are $\lambda$--proper.  Unfortunately, $\lambda$--properness is not in general preserved in iterations; this paper presents a special case where some form of it is.

\begin{defn}\hfill
\begin{enumerate}
\item A set $A\subs P$ is $<\lambda$--linked if every $A_0\in [A]^{<\lambda}$ has a lower bound in $P$.
\item An $(N, P)$--diamond is a sequence $\bar{A}=\langle A_\delta:\delta\in S\rangle$ such that 
\begin{itemize}
\item $S\in \D^+$
\item $A_\delta$ is a subset of $N_\delta\cap P$ with a lower bound in $P$
\item whenever $A\subs N\cap P$ is $<\lambda$--linked, 
\begin{equation}
\{\delta\in S: N_\delta\cap A= A_\delta\}\in\D^+.
\end{equation}
\end{itemize}
\item In the context of (2), if $N_\delta\cap A=A_\delta$ then we say that {\em $\bar{A}$ guesses $A$ at $\delta$.}
\end{enumerate}
\end{defn}

Our first observation is that  that $(N, P)$--diamonds sequences are nothing mysterious -- they are just regular diamond sequences that have been cosmetically altered.

\begin{lem}
Let $N$ be a relevant model with filtration $\langle N_\alpha:\alpha<\lambda\rangle$. Further suppose $\D$ has diamonds.  Then for $S\in\D^+$ we can find an $(N, P)$--diamond $\langle A_\delta:\delta\in S\rangle$.
\end{lem}
\begin{proof}
Let $\langle B_\delta:\delta\in S\rangle$ be a $\D$--diamond sequence, and let $f:\lambda\rightarrow N\cap P$ be a bijection.  Given $\delta\in S$, ask if $f[B_\delta]$ is a $\lambda$--linked subset of $N_\delta\cap P$.  If so, then we let $A_\delta= f[B_\delta]$; if not, then let $A_\delta$ be some arbitrary member of $N\cap P$.

Now suppose $A$ is a $\lambda$--linked subset of $N\cap P$.  Since $\langle B_\delta:\delta\in S\rangle$ is a diamond sequence, we know that the set of $\delta$ for which $B_\delta=f^{-1}(A)\cap\delta$ is in $\D^+$.

There is a closed unbounded set $C$ such that $f\restr\delta$ is a bijection between $\delta$ and $N_\delta\cap P$.  If $\delta\in C$ and $B_\delta= f^{-1}(A)\cap \delta$, then $A_\delta= N_\delta\cap A$.  Since $C\in\D$, we see that $\langle A_\delta:\delta\in S\rangle$ is an $(N, P)$--diamond.
\end{proof}

Starting with the next lemma, we use without mention that the filter $\D$ has a natural interpretation in generic extensions of the universe --- in $V[G]$, we let $\D$ refer to the normal filter generated by $\D\cap V$.

\begin{lem}
Let $\langle A_\delta:\delta\in S\rangle$ be an $(N, P)$--diamond, and let $Q$ be a $\lambda$--complete notion of forcing. If $\dot A$ is a $Q$--name for a $\lambda$--linked subset of $N\cap P$, then 
\begin{equation}
\forces_{Q}\{\delta\in S: N_\delta\cap\dot A= A_\delta\}\in \D^+.
\end{equation}
\end{lem}
\begin{proof}
If not, then we can find a condition $q$ as well as a  $Q$--name $\dot A$ and a sequence $\langle \dot C_i:i<\lambda\rangle$ of $Q$--names such that
\begin{itemize}
\item $\forces_{Q}\dot A\text{ is $\lambda$--linked}$,
\sk
\item $\forces_{Q}\dot C_i\in \D\cap V$, and
\sk
\item $q\forces\delta\in S\cap\bigtriangleup_{i<\lambda}\dot C_i\Longrightarrow A_\delta\neq N_\delta\cap\dot A$.
\end{itemize}

Since $Q$ is $\lambda$--closed, we can find sequences $\langle q_\alpha:\alpha<\lambda\rangle$, $\langle C_\alpha:\alpha<\lambda\rangle$, and $\langle B_\alpha:\alpha<\lambda\rangle$ such that
\begin{itemize}
\item $\alpha<\beta<\lambda\Longrightarrow q_\beta\leq q_\alpha\leq q$ in $Q$
\sk
\item $C_\alpha\in \D$
\sk
\item $B_\alpha\subs N_\alpha\cap P$
\sk
\item $q_\alpha\forces \dot C_\alpha= C_\alpha\text{ and }N_\alpha\cap\dot A=  B_\alpha$
\sk
\end{itemize}
Define $C=\bigtriangleup_{\alpha<\lambda} C_\alpha$.  Since $\D$ is a normal filter, we know that $C\in \D$.

Note that the sequence $\langle B_\alpha:\alpha<\lambda\rangle$ increases with $\alpha$.  Define
\begin{equation}
B=\bigcup_{\alpha<\lambda}B_\alpha.
\end{equation}
It is not hard to see that $B$ is $\lambda$--linked (in the ground model), so there is a $\delta\in S\cap C$ where such that $N_\delta\cap B=A_\delta$. This is a contradiction as $q_\delta$ is an extension of $q$, yet 
\begin{equation}
q_\delta\forces \delta\in S\cap\bigtriangleup_{i<\lambda}\dot C_i\text{ and }N_\delta\cap \dot A=B_\delta = A_\delta.
\end{equation}

\end{proof}
\begin{cor}
If $\bar{A}$ is an $(N, P)$--diamond and $G\subs P$ is a generic subset of $P$, then
\begin{equation}
\{\delta\in S: N_\delta\cap G = A_\delta\}\in \D^+.
\end{equation}
\end{cor}
\begin{proof}
This follows because $G$ is $\lambda$--directed, hence $\lambda$--linked.
\end{proof}

\begin{defn}
A sequence $\bar{R}=\langle (A_\delta, q_\delta):\delta\in S\rangle$ is said to be an $(N, P)$--rule if
\begin{itemize}
\item $\langle A_\delta:\delta\in S\rangle$ is an $(N, P)$--diamond,
\sk
\item $q_\delta$ is a lower bound for $A_\delta$ in $N\cap P$, and 
\sk
\item if $D\in N$ is a dense subset of $P$, then $q_\delta\in D$ for all sufficiently large $\delta\in S$.
\end{itemize}
\end{defn}

\begin{defn}
A notion of forcing $P$ is proper over $\D$--diamonds if for almost every relevant model $N$, whenever we are given an $(N, P)$--rule $\bar{R}=\langle (A_\delta, q_\delta):\delta\in S\rangle$, for every $p\in N\cap P$ there is $q\leq p$ 
\begin{equation*}
q\forces\text{for some  $C\in \D$, if $\delta\in S\cap C$ and $A_\delta\subs \dot G_P$, then $q_\delta\in\dot G_P$.}
\end{equation*}
We say that $q$ is $(N, P, \bar{R})$--generic.
\end{defn}

In other words, $q$ is $(N, P, \bar{R})$--generic if $q$ forces that in the generic extension, for $\D$--almost all $\delta\in S$, if $A_\delta$ guesses $N_\delta\cap G$, then $q_\delta\in G$.  We say that $q$ forces $N\cap G$ to {\em obey the rule $\bar{R}$}. 

\begin{prop}
Suppose $N$ is a relevant model containing $P$, $\bar{R}$ is an $(N, P)$--rule, and $q$ is $(N, P, \bar{R})$--generic.  Then $q$ is $(N, P)$--generic, i.e., 
\begin{equation}
q\forces N[\dot G_P]\cap\ord=N\cap\ord.
\end{equation}
In particular, if $P$ is proper over $\D$--diamonds, then forcing with $P$ preserves the cardinal $\lambda^+$.
\end{prop}

\section{Iterations}

We begin with an outline that shows how properness over $\D$--diamonds is preserved in a simple two--step iteration. Thus, suppose $P$ is proper for $\D$--diamonds and $\forces_{P}\text{``}\dot Q$ is proper for $\D$--diamonds''.  In the following discussion we will show that the composition $P*\dot Q$ is proper for $\D$--diamonds.

Suppose now that $\bar{R}=\langle (A_\delta, p_\delta*\dot q_\delta):\delta\in S\rangle$ is an $(N, P*\dot Q)$--rule and $p\in N\cap P$.  Let $B_\delta$ be the set of ``first co--ordinates'' of members of $A_\delta$.  It is straightforward to prove that $\bar{R}\restr P:=\langle (B_\delta, p_\delta):\delta\in S\rangle$ is an $(N, P)$--rule.  What we want to show is that $(N, P, \bar{R}\restr P)$--generic conditions can be extended in a natural way to $(N, P*\dot Q, \bar{R})$--generic conditions.

To see how this can be accomplished, suppose that $r$ is $(N, P, \bar{R}\restr P)$--generic with $r\leq p$, and $G$ is a generic subset of $P$ containing $r$.  In $V[G]$, let us define
\begin{equation*}
S_0=\{\delta\in S: p_\delta\in G\}.
\end{equation*}
Prior considerations tell us that $S_0\in\D^+$.

Given $\delta\in S_0$, let $B_\delta$ be the set of interpretations of the ``second coordinates'' of members of $A_\delta$, i.e., for $\delta\in S_0$, 
\begin{equation*}
\dot s[G]\in B_\delta\Longleftrightarrow r*\dot s\in A_\delta\text{ for some $r\in P$}.
\end{equation*}
Standard arguments show us that $\bar{R}/G:=\langle (B_\delta, q_\delta[G]):\delta\in S_0\rangle$ is an $(N[G], \dot Q[G])$--rule in $V[G]$.

Now suppose $q\in N[G]\cap\dot Q[G]$ (note that $N[G]\cap\dot Q[G]=N\cap\dot Q[G]$ since $r\in G$).  Since $\dot Q[G]$ is proper for $\D$--diamonds, we can find a condition $s\leq q$ in $\dot Q[G]$ such that $s$ is $(N[G], \dot Q[G], \bar{R}/G)$--generic.  Back in the ground model $V$, we can find a name $\dot s$ forced by $r$ to have the properties ascribed to $s$ in $V[G]$.  It is straightforward to prove that $r*\dot s$ is $(N, P*\dot Q, \bar{R})$--generic and $r*\dot s\leq p*\dot q$.  Thus we have shown that $P*\dot Q$ is proper for $\D$--diamonds.

Now what happens with longer iterations? Assume now that $\mathbb{P}=\langle P_i, \dot Q_i:i<\kappa\rangle$ is $\lambda$--support iteration of $\lambda$--closed notions of forcing such that
\begin{equation}
\forces_{P_i}\dot Q_i\text{ is proper for $\D$--diamonds.}
\end{equation}
We will show that $P_\kappa$, the limit of $\mathbb{P}$, is proper for $\D$--diamonds, so in particular forcing with $P_\kappa$ preserves ~$\lambda^+$.

\begin{thm}[Iteration Theorem]\hfill\\
Let $\langle P_i, \dot Q_i:i<\kappa\rangle$ be a $\lambda$--support iteration such that $\forces_{P_i}\dot Q_i\text{ is proper over $\D$--diamonds.}$ Then $P_\kappa$ is proper over $\D$--diamonds.
\end{thm}

\begin{defn}
Let $N$ be a relevant model with $\mathbb{P}\in N$, and suppose  $i<j$ in $N\cap(\kappa+1)$. Let $\bar{A}=\langle A_\delta:\delta\in S\rangle$ be an $(N, P_j)$--diamond. Given $\delta\in S$, we define
\begin{equation}
A_\delta\restr i = \{p\restr i:p\in A_\delta\}, 
\end{equation}
and 
\begin{equation}
\bar{A}\restr i=\langle A_\delta\restr i:\delta\in S\rangle.
\end{equation}
Similarly, if $\bar{R}=\langle (A_\delta, q_\delta):\delta\in S\rangle$ is an $(N, P)$--rule, we define
\begin{equation}
\bar{R}\restr i = \langle (A_\delta\restr i, q_\delta\restr i):\delta\in S\rangle.
\end{equation}
\end{defn}

\begin{lem}
Let $N$ be a relevant model containing $\mathbb{P}$, and let $i<j$ in $N\cap(\kappa+1)$.
If $\bar{A}$ is an $(N, P_j)$--diamond, then $\bar{A}\restr i$ is an $(N, P_i)$--diamond.
If $\bar{R}$ is an $(N, P_j)$--rule, then $\bar{R}\restr i$ is an $(N, P_i)$--rule.
\end{lem}

\begin{proof}[{\rm\bf Proof of the Iteration Theorem}]\hfill

We prove by induction on $j\in N\cap\kappa+1$ that whenever we are given objects $i$, $\dot p$, and $r$ such that

\begin{itemize}
\item $i<j$
\sk
\item $r\in P_i$
\sk
\item $\forces_{P_i}\dot p\in P_\kappa$
\sk
\item $r\forces\dot p\restr i\in \dot G_{P_i}$
\sk
\item $r$ is $(N, P_i, \bar{R}\restr i)$--generic
\sk
\end{itemize}
we can find a condition $s\in P_j$ such that
\begin{itemize}
\sk
\item $s\restr i= r$
\sk
\item $s$ is $(N, P_j, \bar{R}\restr j)$--generic
\sk
\item $s\forces\dot p\restr j\in\dot G_j$
\sk
\end{itemize}

\bigskip
\noindent{\bf CASE 1}: $j$ is a successor ordinal
\medskip

Let $j=j_0+1$.  Since $j_0$ must be in $N\cap (\kappa+1)$,  we may apply our induction hypothesis to obtain a condition $s_0\in P_{j_0}$ such that
\begin{itemize}
\item $s_0\restr i = r$
\sk
\item $s_0$ is $(N, P_{j_0}, \bar{R}\restr j_0)$--generic, and
\sk
\item $s_0$ forces that $\dot p\restr j_0$ is in $\dot G_{j_0}$.
\sk
\end{itemize}

At this point, we are essentially in the case where we are doing a two--step iteration -- if we view $P_j$ as a two--step iteration $P_{j_0}*\dot Q_{j_0}$, then the arguments presented at the beginning of this section show how to extend $s_0$ to the required $(N, P_j, \bar{R}\restr j)$--generic condition $s$.

\bigskip
\noindent{\bf CASE 2}:  $j$ is a limit ordinal of cofinality $<\lambda$
\medskip

\begin{lem}
\label{lem2}
Suppose $\epsilon\in N\cap (\kappa+1)$ satisfies $\cf(\epsilon)<\lambda$, and we are given sequences $\langle i_\alpha:\alpha<\cf(\epsilon)\rangle$ and $\langle r_\alpha:\alpha<\cf(\epsilon)\rangle$ such that
\begin{itemize}
\item $\langle i_\alpha:\alpha<\cf(\epsilon)\rangle$ is a strictly increasing sequence of ordinals in $N\cap\epsilon$
\item $r_\alpha$ is $(N, P_{i_\alpha}, \bar{R}\restr i_\alpha)$--generic
\item $\alpha<\beta<\kappa\Longrightarrow r_\beta\restr i_\alpha= r_\alpha$.
\end{itemize}
Then the condition $s:=\bigcup_{\alpha<\cf(\epsilon)} r_\alpha$ is $(N, P_\epsilon, \bar{R}\restr \epsilon)$--generic.
\end{lem}
\begin{proof}
Clearly $s\in P_\epsilon$ as we are using $\lambda$--support iteration. Let $G$ be any generic subset of $P_\epsilon$ that contains $s$; we will work in the generic extension $V[G]$.

For $\alpha<\cf(\epsilon)$, let $G_\alpha= G\restr P_{i_\alpha}$.  Clearly $r_\alpha\in G_\alpha$ and $G_\alpha$ is a generic subset of $P_{i_\alpha}$, so there is a set $C_\alpha\in \D$ such that
\begin{equation}
\label{eqn1}
\delta\in S\cap C_\alpha\text{ and }N_\delta\cap G_\alpha= A_\alpha\restr i_\alpha\Longrightarrow q_\delta\restr i_\alpha\in G_\alpha.
\end{equation}
Let $C=\bigcap_{\alpha<\cf(\epsilon)}C_\alpha\in \D$.  Given $\delta\in S\cap C$ if $\bar{A}$ guesses $G$ at $\delta$, then (\ref{eqn1}) implies that $q_\delta\restr i_\alpha\in G_\alpha$ for all $\alpha<\cf(\epsilon)$.  Since $G$ is a generic subset of $P_\epsilon$, it follows that $q_\delta$ is in $G$, as required.
\end{proof}

Now we return to the case where $\cf(j)<\lambda$.  Let $\langle i_\alpha:\alpha<\cf(j)\rangle$ be increasing, continuous, and cofinal in $N\cap j$ --- note that we can achieve continuity because $N$ is closed under sequences of length $<\lambda$. Without loss of generality we assume $i_0=i$.  

By induction on $\alpha<\cf(j)$, we choose  conditions $r_\alpha\in P_{i_\alpha}$ such that

\begin{itemize}
\item $r_0=r$
\sk
\item $r_\alpha\forces\dot p_\alpha\restr i_\alpha\in \dot G_{P_{i_\alpha}}$
\sk
\item if $\beta<\alpha$ then $r_\alpha\restr i_{\beta}= r_\beta$
\sk
\item if $\alpha$ is a limit ordinal, then $r_\alpha=\bigcup_{\beta<\alpha} r_\beta$
\sk
\item $r_\alpha$ is $(N, P_{i_\alpha},\bar{R}\restr i_\alpha)$--generic
\sk
\end{itemize}

The construction of $\langle r_\alpha:\alpha<\cf(j)\rangle$ is straightforward --- at successor stages we apply our induction hypothesis, while at limit stages we invoke Lemma \ref{lem2} to show that the construction continues.

Another application of Lemma \ref{lem2} shows us that $s$ is $(N, P_j, \bar{R}\restr j)$--generic; the other requirements for $s$ are also easily verified.

\bigskip
\noindent{\bf CASE 3}: $\cf(j)=\lambda$
 \medskip

Let $\langle i_\alpha:\alpha<\lambda\rangle$ be increasing, continuous, and cofinal in $N\cap j$ with $i_0=i$.  Let $\langle D_\alpha:\alpha<\lambda\rangle$ list all dense open subsets of $P_j$ that are elements of $N$.

By induction on $\alpha<\lambda$, we will define objects $\dot p_\alpha$ and $r_\alpha$ such that

\begin{enumerate}
\item $r_0=r$, $\dot p_0=\dot p\restr j$
\sk
\item $r_\alpha$ is $(N, P_{i_\alpha},\bar{R}\restr i_\alpha)$--generic
\sk
\item $r_\alpha\restr i_\beta= r_\beta$ for $\beta<\alpha$
\sk
\item $r_\alpha\forces\dot p_\alpha\in N\cap P_j\text{ and } \dot p_\alpha\restr i_\alpha\in \dot G_{P_{i_\alpha}}$
\sk
\item $r_{\alpha+1}\forces \dot p_{\alpha+1}\in D_\alpha$
\sk
\item for $\beta<\alpha$, $r_\alpha\forces\dot p_\alpha\leq\dot p_\beta$
\sk
\item for $\alpha\in S$, $r_\alpha$ forces the statement
\begin{equation*}
(\otimes)\quad \text{if } q_\alpha\restr i_\alpha\in\dot G_{i_\alpha}\text{ and $q_\alpha$ is a lower bound for $\langle \dot p_\beta:\beta<\alpha\rangle$, then $\dot p_\alpha=q_\alpha\restr j$.}
\end{equation*}
\sk
\sk
\end{enumerate}

\sk
\sk

\noindent{\bf Construction of $\langle \dot p_\alpha:\alpha<\lambda\rangle$ and $\langle r_\alpha:\alpha<\lambda\rangle$}:
\sk
\sk

\noindent{\bf Initial stage}:
\sk
We have already defined $r_0$ and $\dot p_0$.
\sk
\sk
\noindent{\bf Successor stages}:
\sk

Assume now that $\alpha$ is a successor ordinal, say $\alpha=\beta+1$.
Our construction will give us objects $r_\beta$ and $\dot p_\beta$ satisfying the appropriate conditions.  We apply our induction hypothesis with $i_\alpha$, $i_\beta$, $\dot p_\beta\restr i_\alpha$, $r_\beta$, and $\bar{R}\restr i_\alpha$ standing for the objects $j$, $i$, $\dot p$, $r$, and $\bar{R}$ appearing there.   This gives us an object $r_\alpha$ such that
\begin{itemize}
\sk
\item $r_\alpha$ is $(N, P_{i_\alpha}, \bar{R}\restr i_\alpha)$--generic,
\sk
\item $r_\alpha\restr i_\beta= r_\beta$, and
\sk
\item $r_\alpha\forces\dot p_\beta\restr i_\alpha\in\dot G_{i_\alpha}$. 
\sk
\end{itemize}

Now let $G$ be any generic subset of $P_{i_\alpha}$ that contains $r_\alpha$. We know that $N\cap G$ is $P_{i_\alpha}$--generic over $N$ because $r_\alpha$ is $(N, P_{i_\alpha})$--generic. Since $D_\beta\in N$, a standard genericity argument tells us that there is a condition $p_{\alpha}\in N[G]\cap P_j= N\cap P_j$ such that
\begin{itemize}
\item $p_{\alpha}\restr i_\alpha\in G$,
\sk
\item $p_{\alpha}\leq \dot p_\beta[G]$, and
\sk
\item $p_{\alpha}\in D_\beta$.
\sk
\end{itemize}
Back in $V$, we let $\dot p_\alpha$ be a name for this $p_\alpha$; it should be clear that $\dot p_\alpha$ is as required.
\sk
\sk
\noindent{\bf Limit stages}:
\sk

If $\alpha$ is a limit ordinal, we know
\begin{equation*}
r_\alpha=\bigcup_{\beta<\alpha} r_\beta.
\end{equation*}
Since $\cf(\alpha)<\lambda$, Lemma \ref{lem2} implies that $r_\alpha$ is $(N, P_{i_\alpha}, \bar{R}\restr i_\alpha)$--generic.  Also, our inductive assumptions imply that for all $\beta<\alpha$,
\begin{equation*}
r_\alpha\forces\dot p_\beta\restr i_\alpha\in\dot G_{i_\alpha}.
\end{equation*}

Let $G$ be any generic subset of $P_{i_\alpha}$ with $r_\alpha\in G$.  In the extension $V[G]$, each name $\dot p_\beta$ is interpreted as a condition in $N\cap P_j$, and we know
\begin{itemize}
\item $\forall\beta<\alpha$, $p_\beta\restr i_\alpha\in G$, and
\sk
\item $\langle p_\beta:\beta<\alpha\rangle$ is decreasing.
\end{itemize}
Now we ask the question
\sk
\sk
\noindent{\em Is it the case that
\begin{itemize}
\item $\alpha\in S$
\item $q_\alpha\restr i_\alpha\in G$, and
\item $q_\alpha\restr j$ is a lower bound for $\langle p_\beta:\beta<\alpha\rangle$ in $P_j$?
\end{itemize}
}

\sk
\sk

If the answer is yes, then we let $p_\alpha=q_\alpha\restr j$.  If the answer is no, then we let $p_\alpha$ be a lower bound for $\langle p_\beta:\beta<\alpha\rangle$ in $N\cap P_j$ with $p_\alpha\restr i_\alpha\in G$.

Now back in the ground model, we let $\dot p_\alpha$ be a name forced by $r_\alpha$ to be as above.  Note that $\dot p_\alpha$ is as required in $(\otimes)$, and our construction continues.
\bigskip

Once we have defined $r_\alpha$ and $\dot p_\alpha$ for every $\alpha<\lambda$, we define
\begin{equation*}
s:=\bigcup_{\alpha<\lambda}r_\alpha.
\end{equation*}
Clearly $s\restr i = r$ and $s\forces \dot p\restr j \in\dot G_j$, so we need only verify that $s$ is $(N, P_j, \bar{R}\restr j)$--generic.

Let $G$ be any generic subset of $P_j$ that contains $s$, and step into the model $V[G]$. Each $\dot p_\alpha$ is interpreted as some $p_\beta\in N\cap P_j$ and our construction guarantees that the filter generated by $\langle p_\alpha:\alpha<\lambda\rangle$ is generic over $N$ and hence equal to $N\cap G$.
This tells us that $s$ is $(N, P_j)$--generic.

For each $\alpha<\lambda$, the condition $r_\alpha$ is $(N, P_{i_\alpha}, \bar{R}\restr i_\alpha)$--generic so in $V[G]$ we can find a set $C_\alpha\in \D$ that witnesses this, i.e., if $\delta\in C_\alpha\cap S$ and $q_\delta\restr i_\alpha$ guesses $N_\delta\cap G\restr i_\alpha$, then $q_\delta\restr i_\alpha\in G\restr i_\alpha$.

Since $\langle p_\alpha:\alpha<\lambda\rangle$ generates $N\cap G$ and $N\cap G$ is generic over $N$, there is a closed unbounded set $E\subs\lambda$ such that
\begin{equation}
\delta\in E\Longrightarrow \langle p_\alpha:\alpha<\delta\rangle\text{ generates a generic subset of $N_\alpha\cap P$}.
\end{equation}
Let $C=E\cap\bigtriangleup_{\alpha<\lambda} C_\alpha$; since $\D$ is normal, we know that $C\in \D$.

\begin{claim}
If $\delta\in C\cap S$ and $q_\delta\restr j$ guesses $N_\delta\cap G$, then $q_\delta\restr j\in G$.
\end{claim}
\begin{proof}
Suppose we are given such a $\delta$. It suffices to show that $q_\delta\restr i_\delta\in G_{i_\delta}$ and $q_\delta\restr j$ is a lower bound for $\langle p_\beta:\beta<\delta\rangle$ --- if this happens, then our construction guarantees $p_\delta=q_\delta\restr j$ and $p_\delta\in G$.

Our definition of $C$ implies that $\delta\in C_\beta$ for all $\beta<\delta$. Since $q_\delta\restr j$ guesses $N_\delta\cap G$, we know that $q_\delta\restr i_\alpha$ guesses $N\cap G_{i_\alpha}$ for all $\alpha<\lambda$.  Given $\beta<\delta$, we know that $r_\beta\in G_{i_\beta}$ and $r_\beta$ is $(N, P_{i_\beta}, \bar{R}\restr i_\beta)$--generic.  Putting all this together, we may conclude that for all $\beta<\delta$, $q_\delta\restr i_\beta\in G_{i_\beta}$, hence $q_\delta\restr i_\delta \in G_{i_\delta}$.

Now why is $q_\delta\restr j$ a lower bound for $\langle p_\beta:\beta<\delta\rangle$?  This follows because $\delta\in C$ --- the sequence $\langle p_\beta:\beta<\delta\rangle$ generates $N_\delta\cap G$, and we have assumed that $q_\delta\restr j$ guesses $N_\delta\cap G$.

Since $r_\delta$ forces $(\otimes)$ to hold, we know that $\dot p_\delta[G]=q_\delta\restr j$, hence $q_\delta\restr j\in G$.

\end{proof}

We have therefore shown that $s$ is $(N, P_j, \bar{R}\restr j)$--generic. Since $s\restr i= r$ and our construction guarantees that $s\forces\dot p\restr j\in \dot G_{P_j}$, so $s$ is as required.

\bigskip

\noindent{\bf CASE 4}: $\cf(j)>\lambda$

\medskip

The construction in this case is very similar to that of the previous case.   Let $k=\sup(N\cap j)$; since $N$ is closed under sequences of length $<\lambda$, it follows that $\cf(k)=\lambda$ and we can fix a continuous increasing sequence $\langle i_\alpha:\alpha<\lambda\rangle$ of elements of $N\cap j$ cofinal in $k$.

The idea now is to mimic the construction given for the case where $\cf(j)=\lambda$.  Let $\langle D_\alpha:\alpha<\lambda\rangle$ list all dense open subsets of $P_j$ that are elements of $N$.  By induction on $\alpha<\lambda$, define objects $\dot p_\alpha$ and $r_\alpha$ satisfying exactly the same requirements as in the previous case --- that construction did not require that $j$ was an element of $N\cap \kappa$, only that a sequence along the lines of $\langle i_\alpha:\alpha<\lambda\rangle$ exists.  One then checks that the resulting condition $s$ defined as there has all the required properties.  Note that what's going on is that members of $N\cap P_j$ are actually members of $N\cap P_k$ --- the support of a condition in $N\cap P_\kappa$ is a subset of $N\cap\kappa$ because $\lambda\subs N$.

\end{proof}

\section{An Example}

Let $S\subs S^{\omega_2}_{\omega_1}:=\{\delta<\omega_2:\cf(\delta)=\omega_1\}$ be stationary.  Recall that a continuous ladder system on $S$ is a family of functions $\bar{\eta}=\langle \eta_\delta:\delta\in S\rangle$ such that $\eta_\delta$ is a strictly increasing and continuous from $\omega_1$ onto a cofinal subset of $\delta$.

A continuous ladder system $\bar{\eta}$ has the club uniformization property if whenever $\bar{c}=\langle c_\delta:\delta\in S\rangle$ is a family of functions from $\omega_1$ to $\{0,1\}$, there is a function $h$ such that for all $\delta\in S$, the set $\{i<\omega_1\: c_\delta(i)= h_\delta(i)\}$ contains a closed unbounded subset of $\omega_1$.

Shelah \cite{Sh:f} has shown that if the Continuum Hypothesis is true, then no continuous ladder system on (all of) $S^{\omega_2}_{\omega_1}$ has the club uniformization property.  If we are looking at a stationary $S\subs  S^{\omega_2}_{\omega_1}$ such that $S^{\omega_2}_{\omega_1}\setminus S$ is stationary as well, then the techniques of \cite{Sh:587} show how to build a model where the Continuum Hypothesis holds and continuous ladder systems on $S$ have the club uniformization property. 

Let us fix a stationary, co--stationary $E_0\subs\omega_1$ and let $\D$ be the club filter restricted to $\omega_1\setminus E_0$.  Further assume that $\mathfrak{D}$ has diamonds --- this follows if $V=L$ or if, e.g.,  $\diamondsuit^*(\omega_1\setminus E_0)$ holds.  

We will force a weak version of the club uniformization property to hold for a continuous ladder system $\bar{\eta}=\langle\eta_\delta:\delta\in S\rangle$ on $S:=S^{\omega_2}_{\omega_1}$; what we achieve is that for every family $\bar{c}=\langle c_\delta:\delta\in S\rangle$ of functions mapping $\omega_1$ to $\{0, 1\}$, there is a function $h:\omega_2\rightarrow 2$ such that for each $\delta\in S$, 
\begin{equation}
\label{eqn2}
\{i\in E_0: h(\eta_\delta(i))\neq c_\delta(i)\}\text{ is non--stationary.}
\end{equation}
Said another way, for each $\delta\in S$ there is a closed unbounded $C_\delta\subs\omega_1$ such that
\begin{equation}
i\in C_\delta\cap E_0\Longrightarrow h(\eta_\delta(i))=c_\delta(i);
\end{equation}
i.e., $h$ achieves success at almost every point in $\eta_\delta[E_0]$.

Let us fix a continuous ladder system $\bar{\eta}=\langle \eta_\delta:\delta\in S\rangle$.  Suppose $\langle c_\delta:\delta\in S\rangle$ is a family of functions each mapping $\omega_1$ to $\{0, 1\}$.  Our first goal is to define a notion of forcing that will adjoin a function $h$ such that (\ref{eqn2}) is satisfied for all $\delta\in S$.

A condition $p$ is simply an approximation to the desired $h$ of size $\leq\omega_1$, i.e.,  $p\in P$ if $p$ is a function satisfying
\begin{itemize}
\sk
\item $\dom(p)\in [\omega_2]^{\leq\omega_1}$
\sk
\item $\ran(p)\subs\{0, 1\}$
\sk
\item for all $\delta\in S$, 
\begin{equation*}
\{i\in E_0:  p(\eta_\delta(i))\neq c_\delta(i)\}\text{ is non--stationary.}
\end{equation*}
\sk
\end{itemize}

Clearly $P$ is $<\omega_1$--closed and for each $\alpha<\omega_2$, the set of conditions with $\alpha$ in their domain is dense in $P$.  Thus forcing with $P$ adds no new countable sequences to the ground model and adjoins a function from $\omega_2$ to $\{0, 1\}$. 

\begin{claim}
$P$ is proper for $\D$--diamonds.
\end{claim}
\begin{proof}

Let $N$ be a relevant model with filtration $\langle N_i:i<\omega_1\rangle$ and let $p\in N\cap P$ be arbitrary.  Suppose $E_1\in\D^+$ and let $\bar{R}=\langle (A_\delta, q_\delta):\delta\in S\rangle$ be an $(N, P)$--rule. Note that we may assume that $E_0\cap E_1=\emptyset$ because of our definition of $\D$. We will construct a decreasing sequence $\langle p_\alpha:\alpha<\omega_1\rangle$ of conditions in $N\cap P$ in such a way that $q:=\bigcup_{\alpha<\omega_1}p_\alpha$ is an $(N, P, \bar{R})$--generic extension of $p$.

Let $\gamma=N\cap\omega_2$, and for $\alpha<\omega_1$ let $\gamma_\alpha=N_\alpha\cap\omega_2$.  The sequence $\langle\gamma_\alpha:\alpha<\omega_1\rangle$ is strictly increasing, continuous, and cofinal in $\gamma$.

As we build the sequence $\langle p_\alpha:\alpha<\omega_1\rangle$, we will also be defining a strictly increasing and continuous sequence of countable ordinals $\langle i_\alpha:\alpha<\omega_1\rangle$.

We begin by letting $i_0$ be the least $i<\omega_1$ such that $p\in N_i$, and let $p_0\in N\cap P$ be some totally $(N_{i_0}, P)$--generic extension of $p$.

Given $\langle p_\beta:\beta\leq\alpha\rangle$ and $\langle i_\beta:\beta\leq \alpha\rangle$, we let $i_{\alpha+1}$ be the least ordinal $i$ such that both $\langle p_\beta:\beta\leq\alpha\rangle$ and $\langle i_\beta:\beta\leq \alpha\rangle$ are elements of $N_i$.  Note that such an $i$ exists because $N^{<\omega_1}\subs N$.  We let $p_{\alpha+1}$ be a totally $(N_{i_{\alpha+1}}, P)$--generic extension of $p_\alpha$ in $N\cap P$.

Now what happens at limit stages of the construction?  If $\alpha$ is a limit ordinal, we will be handed $\langle p_\beta:\beta<\alpha\rangle$ and $\langle i_\beta:\beta<\alpha\rangle$.  We are committed to the continuity of $\langle i_\alpha:\alpha<\omega_1\rangle$, so this means that we are forced to choose
\begin{equation*}
i_\alpha=\bigcup_{\beta<\alpha} i_\beta.
\end{equation*}
Let us define
\begin{equation*}
r_\alpha=\bigcup_{\beta<\alpha}p_\beta.
\end{equation*}
Since $\alpha$ is a countable ordinal, we know that $r_\alpha$ is a condition in $P$, and the relevance of the model $N$ implies that $r_\alpha\in N\cap P$.  By our construction, we know that $r_\alpha$ is totally $(N_{i_\alpha}, P)$--generic --- this follows because
\begin{equation*}
N_{i_\alpha}=\bigcup_{\beta<\alpha} N_{i_\beta}.
\end{equation*}

\noindent Now we ask:

 \noindent{\em Is it the case that
\begin{itemize}
\sk 
\item $i_\alpha=\alpha$,
\sk
\item $\gamma_\alpha = \eta_\gamma(i_\alpha)$, and  
\sk
\item $\alpha\in E_0\cup E_1$?
\sk
\end{itemize}
}

If not, we let $p_\alpha=r_\alpha$ and the construction continues. If the answer is yes, then we have two cases to consider --- the case $\alpha\in E_0$ and the case $\alpha\in E_1$

If $\alpha\in E_0$, we note first that $\dom(r_\alpha)\subs \gamma_\alpha$ --- this is because $p_\beta\in N_\alpha$ for all $\beta<\alpha$ and $\dom(r_\alpha)=\cup_{\beta<\alpha}\dom(p_\beta)$. Thus we may define
\begin{equation*}
p_\alpha=r_\alpha\cup\{\langle \delta_\alpha, c_\gamma(\alpha)\rangle\},
\end{equation*}
and conclude that $p_\alpha\in N\cap P$. 

If $\alpha\in E_1$, then we ask if $A_\alpha$ is equal to the filter on $N_\alpha\cap P$ generated by $\langle p_\beta:\beta<\alpha\rangle$.  If yes, then we let $p_\alpha=q_\alpha$  (note that $q_\alpha\leq r_\alpha$ if this happens); if not, we let $p_\alpha= r_\alpha$.

In either case, the condition $p_\alpha$ will be in $N\cap P$ and the construction can continue.

\begin{claim}
The sequence $\langle p_\alpha:\alpha<\omega_1\rangle$ has a lower bound in $P$.
\end{claim}
\begin{proof}
Let $q=\bigcup_{\alpha<\omega_1} p_\alpha$.  It is clear that $q$ is a partial function from $\omega_2$ to $\{0, 1\}$ with domain a set of cardinality $\aleph_1$.   Since each $p_\alpha$ is an element of $N$, we know that $\dom(q)\subs\gamma$. 

What we need to show is that for every $\delta\in S$, (\ref{eqn2}) holds. If $\delta>\gamma$, then (\ref{eqn2}) holds because $\dom(q)\subs\gamma$.  If $\delta<\gamma$, we note that $\delta\in N$ (as $N^{<\omega_1}\subs N$ implies $N\cap\omega_2$ is an initial segment of $\omega_2$), and the set of conditions whose domain includes $\delta\cup\{\eta_\delta(i):i<\omega_1\}$ is dense in $P$ and an element of $N$.  Thus there is a stage $\alpha$ such that
\begin{equation*}
\delta\cup\{\eta_\delta(i):i<\omega_1\}\subs\dom(p_\alpha).
\end{equation*}
Since $p_\alpha\in P$, the definition of $q$ implies (\ref{eqn2}) holds for $\delta$.

The last case to consider is when $\delta=\gamma$.  Note that there is a closed unbounded set of $\alpha<\omega_1$ for which $i_\alpha=\alpha$ and $\eta_{\gamma}(\alpha)=\gamma_\alpha$.  If $\alpha\in E_0$ has these properties, then at stage $\alpha$ we ensured that $q(\eta_\gamma(\alpha))=c_\gamma(\alpha)$.  Thus (\ref{eqn2}) holds for $\gamma=\delta$, and we have established that $q$ is a condition in $P$.
\end{proof}

\begin{claim}
The condition $q$ is $(N, P, \bar{R})$--generic.
\end{claim}
\begin{proof}
Again, there is a closed unbounded set of $\alpha$ for which $i_\alpha=\alpha$ and $\eta_{\gamma}(\alpha)=\gamma_\alpha$.  Note that for such an $\alpha$, we automatically achieve that $\langle p_\beta:\beta<\alpha\rangle$ generates an $(N_\alpha, P)$--generic filter $G_\alpha$ --- this follows because $N_\alpha=\bigcup_{\beta<\alpha}N_{i_\beta}$.  If for such an $\alpha$ it happens that $G_\alpha= A_\alpha$, then we made sure that $p_\alpha= q_\alpha$.  Since
\begin{equation*}
q\forces N\cap\dot G_{P}\text{ is generated by }\langle p_\alpha:\alpha<\omega_1\rangle,
\end{equation*}
we have ensured that $q$ is $(N, P, \bar{R})$--generic.
\end{proof}
\end{proof}

\def\germ{\frak} \def\scr{\cal} \ifx\documentclass\undefinedcs
  \def\bf{\fam\bffam\tenbf}\def\rm{\fam0\tenrm}\fi 
  \def\defaultdefine#1#2{\expandafter\ifx\csname#1\endcsname\relax
  \expandafter\def\csname#1\endcsname{#2}\fi} \defaultdefine{Bbb}{\bf}
  \defaultdefine{frak}{\bf} \defaultdefine{mathfrak}{\frak}
  \defaultdefine{mathbb}{\bf} \defaultdefine{mathcal}{\cal}
  \defaultdefine{beth}{BETH}\defaultdefine{cal}{\bf} \def\bbfI{{\Bbb I}}
  \def\mbox{\hbox} \def\text{\hbox} \def\om{\omega} \def\Cal#1{{\bf #1}}
  \def\pcf{pcf} \defaultdefine{cf}{cf} \defaultdefine{reals}{{\Bbb R}}
  \defaultdefine{real}{{\Bbb R}} \def\restriction{{|}} \def\club{CLUB}
  \def\w{\omega} \def\exist{\exists} \def\se{{\germ se}} \def\bb{{\bf b}}
  \def\equivalence{\equiv} \let\lt< \let\gt>
\providecommand{\bysame}{\leavevmode\hbox to3em{\hrulefill}\thinspace}
\providecommand{\MR}{\relax\ifhmode\unskip\space\fi MR }
\providecommand{\MRhref}[2]{%
  \href{http://www.ams.org/mathscinet-getitem?mr=#1}{#2}
}
\providecommand{\href}[2]{#2}

\end{document}